\documentclass[a4paper,10pt]{amsart}
\usepackage[utf8]{inputenc}
\usepackage{lmodern}
\usepackage[english]{babel}
\usepackage[T1]{fontenc}
\usepackage{microtype} 
    \usepackage[OT2,OT1]{fontenc}
    \newcommand\cyr{%
    \renewcommand\rmdefault{wncyr}%
    \renewcommand\sfdefault{wncyss}%
    \renewcommand\encodingdefault{OT2}%
    \normalfont
    \selectfont}
    \DeclareTextFontCommand{\textcyr}{\cyr}
\usepackage{amsmath,amssymb,amsfonts,amsthm,amscd,latexsym,mathrsfs}
\usepackage[all]{xy}
\input xypic
\usepackage{color}
\usepackage{hyperref}
\usepackage{lipsum}
\usepackage{breqn}
\hypersetup{colorlinks=true,linkcolor=red,citecolor=blue}
\usepackage{booktabs}
\usepackage{multirow}
\theoremstyle{plain}
\newtheorem{theorem}[subsection]{{\bf Theorem}}
\newtheorem*{theorem*}{{\bf Theorem}}

\newtheorem*{corollary*}{{\bf Corollary}}
\newtheorem{proposition}[subsection]{{\bf Proposition}}
\newtheorem{lemma}[subsection]{{\bf Lemma}}

\theoremstyle{definition}

\theoremstyle{remark}

\numberwithin{equation}{section}
\newcommand{\gen}[1]{\langle #1 \rangle}

\begin{document}
\baselineskip=14pt
\title{On the order of the Schur multiplier of $p$-groups}
%%%%%%%%%%%%%%%%%%%%%%%%%%%%%%%%%%%%%%%%%%%%%%%%%%%%%
%\author[U. Jezernik]{Urban Jezernik}
%\address[Urban Jezernik]{
%Institute of Mathematics, Physics, and Mechanics \\
%Jadranska 19 \\
%1000 Ljubljana \\
%Slovenia}
%\email{urban.jezernik@imfm.si}
%%%%%%%%%%%%%%%%%%%%%%%%%%%%%%%%%%%%%%%%%%%%%%%%%%%%%
\author[P. K. Rai]{Pradeep K. Rai}
\address[Pradeep K. Rai]{Department of Mathematics, Bar-Ilan University
Ramat Gan \\
Israel}
\email{raipradeepiitb@gmail.com}
%%%%%%%%%%%%%%%%%%%%%%%%%%%%%%%%%%%%%%%%%%%%%%%%%%%%%
%%%%%%%%%%%%%%%%%%%%%%%%%%%%%%%%%%%%%%%%%%%%%%%%%%%%%
\begin{abstract}
We give a bound on the order of the Schur multiplier of $p$-groups refining earlier bounds. As an application we complete the classification of groups having Schur multiplier of maximum order. Finally we prove that the order of the Schur multiplier of a finite $p$-group of maximal class and order $p^n$ is at most $p^{\frac{n}{2}}$.
%Let $G$ be a finite $p$-group of order $p^n$ and maximal class and $M(G)$ denote the Schur multiplier of $G$. We prove that $|M(G)| \leq p^{\frac{n}{2}}$.
%A classical result of Green states that $|M(G)| \leq p^{\frac{1}{2}n(n-1)}$. In 2009, Niroomand, improving Green's and other bounds on $|M(G)|$ for a non-abelian $p$-group $G$, proved that $|M(G)| \leq p^{\frac{1}{2}(n-k-1)(n+k-2)+1}$. In this article we prove that a bound, obtained earlier, by Ellis and Wiegold is more general than the bound of Niroomand. We derive from the bound of Ellis and Wiegold that $|M(G)| \leq p^{\frac{1}{2}(d(G)-1)(n+k-2)+1}$ for a non-abelian $p$-group $G$. Moreover, we sharpen the bound of Ellis and Wiegold and as a consequence derive that if $G^{ab}$ %\neq C_{p^e} \times C_{p^e} \times \ldots  \times C_{p^e}$, 
%is not homocyclic then $|M(G)| \leq p^{\frac{1}{2}(d(G)-1)(n+k-3)+1}$. We further note an improvement in an old bound given by Vermani. Finally we prove, for a $p$-group of coclass $r$,  that $|M(G)| \leq p^{\frac{1}{2}(r^2-r)+kr+1}$. This improves a bound by Moravec.  %Finally we consider the corank $t$ of a group $G$, defined by $|M(G)| = p^{\frac{1}{2}n(n-1)-t}$. The groups with small values of $t$ has been classified by many authors. We note a bound on $n$ in terms of $t$, in particular we note that $n \leq t+2$. % prove for a finite $p$-group $G$ of coclass $r$, $|M(G)| \leq p^{(r^2-r)/2+kr-h_G}$. This significantly improves a bound of Moravec \cite{Moravec1}. We also improve a bound on $|M(G)$ given by Vermani \cite{Vermani2}, if the nilpotency class of $G$ is at least 3.
\end{abstract}
%%%%%%%%%%%%%%%%%%%%%%%%%%%%%%%%%%%%%%%%%%%%%%%%%%%%%
\maketitle
%%%%%%%%%%%%%%%%%%%%%%%%%%%%%%%%%%%%%%%%%%%%%%%%%%%%%
%%%%%%%%%%%%%%%%%%%%%%%%%%%%%%%%%%%%%%%%%%%%%%%%%%%%%
\section{Introduction}

Let $G$  be a group. The center and the commutator subgroup of $G$ are denoted by $Z(G)$, and $\gamma_2(G)$ respectively. By $d(G)$ we denote the minimal no of generators of $G$. We write $\gamma_i(G)$ and $Z_i(G)$ for the $i$-th term in the lower and upper central series of $G$ respectively. Finally, the abelianization of the group $G$, i.e. $G/\gamma_2(G)$, is denoted by $G^{ab}$.

 %and $x,y \in G$. Then $[x,y]$ denotes the commutator $x^{-1}y^{-1}xy$. By $\gen{x}$ we denote the cyclic subgroup generated by $x$. Let $H$ be a subset of  $G$. Then $\gen{H}$ denotes the subgroup generated by $H$. 
 %The group of all homomorphisms from a group $H$ to a group $K$ is denoted by $Hom(H, K)$. Let $f \in Hom(H,K)$, then $ker(f)$ denotes the kernal of $f$. 

%Let $p$ be a prime number. By $C_{p^n}$ we denote the cyclic group of order $p^n$, and by $C_{p^n}^{(k)}$, the direct product of $k$ coppies of $C_{p^n}$. The subgroup $\gen{x^{p} \mid x \in G}$ is denoted by $G^{p}$. 
%If $H$ is a subgroup (proper subgroup) of $G$, then we write $H \le G$ ($H < G$). %If $H$ is a subset (proper subset) of $G$, then we write $H \subseteq G$ ($H \subset G$). 

%The subgroup of $G$ generated by all the elements of order $p$ is denoted by $\Omega_1(Z(G))$
%, $x^a$ denotes the conjugate of $x$ by $a$, i.e., $a^{-1}xa$, $x^G$ denotes the $G$-conjugacy class of $x$ and $[x, G]$ denotes the set of all $[x,z]$, for $z \in G$. %and $G_{p^i} = \gen{x \in G \mid x^{p^i} = 1}$, where $i \ge 1$ is an integer. 

\vspace{.2cm}
Let $G$ be finite $p$-group of order $p^n$ and $M(G)$ be the Schur multiplier of $G$. In 1956 Green proved that $|M(G)| \leq p^{\frac{1}{2}n(n-1)}$ \cite{Green}. Since then, Green's bound has been refined by many mathematicians. For a detailed account %of the results on the bound on $M(G)$ 
we refer the reader to a note recently written by the author \cite{Rai1}. %This earlier note misses to note a comment in a paper written by Ellis and Weigold \cite{Ellis1}. It turns out that this comment can play a crucial role when dealing with the groups with higher nilpotency classes. 
There has also been interest in the  classification of finite $p$-groups given the order of their Schur multiplier. By Green's result note that $|M(G)| = p^{\frac{1}{2}n(n-1)-t}$ for some non-negative integer $t$. The groups with $t= 0, 1, 2 ,3, 4, 5$ have been classified, see \cite{Berkovich, Ellis2,  Salemkar, Zhou, Niroomand4}. The first such result was proved by Berkovich \cite{Berkovich}. He proved that $|M(G)| = p^{\frac{1}{2}n(n-1)}$ if and only if $G$ is elementary abelian.

Let G be a non-abelian $p$-group of order $p^n$ with $|\gamma_2(G)| = p^k$ and $d(G) = d$. Niroomand proved in \cite{Niroomand2} that 
\begin{equation} \label{eq0}
|M(G)| \leq p^{\frac{1}{2}(n-k-1)(n+k-2)+1}. 
\end{equation}

The author noted in \cite{Rai1} that a bound of Ellis and Weigold is actually better than this bound and derived from their bound that 
\[|M(G)| \leq p^{\frac{1}{2}(d-1)(n+k-2)+1}(= p^{\frac{1}{2}(d-1)(n+k)-(d-2)}).\] %This is further refined in the following

Niroomand also classified finite $p$-groups such that $k=1$ and bound (\ref{eq0}) is attained. Note that if $k = 1$, then the group is of nilpotency class 2. The author further classified the groups of nilpotency class 2 such that this bound is attained \cite{Rai2}. Recently Hatui proved that there are no $p$-groups, for $p \geq 5$, of nilpotency class 3 or more attaining the bound \cite{Hatui}. She also gave an example of 3-group of nilpotency class 3 such that the bound is attained. 

A natural question is: whether this bound can be improved for groups of higher nilpotency classes. We offer the following theorem which improves the bound for the groups with large $\min(d(G),c)$, where $c$ is the nilpotency class of $G$.

\begin{theorem}\label{thm1}
Let G be a non-abelian $p$-group of order $p^n$ and nilpotency class $c$ with $|\gamma_2(G)| = p^k$ and $d(G) = d$. Then
$$|M(G)| \leq p^{\frac{1}{2}(d-1)(n+k)-\mathop{\sum_{i =2}^{\min(d,c)} d-i}}.$$ 
\end{theorem}

%The interest in noting a bound containing the term nilpotency class arises due to the interest of many authors in the  classification of finite $p$-groups given the order of their Schur multiplier. By Green's result note that $|M(G)| = p^{\frac{1}{2}n(n-1)-t}$ for some nonnegative integer $t$. The groups with $t= 0, 1, 2 ,3, 4, 5$ have been classified, see \cite{Berkovich, Ellis2,  Salemkar, Zhou, Niroomand4}. The first such result was proved by Berkovich \cite{Berkovich}. He proved that $|M(G)| = p^{\frac{1}{2}n(n-1)}$ if and only if $G$ is elementary abelian.  

%Niroomand classified finite $p$-groups such that $k=1$ and the Bound \ref{eq0} is attained. The author further classified the groups of nilpotency class 2 such that this Bound is attained \cite{Rai2}. Recently Hatui proved that There are no $p$-groups ($p \neq 3$) of nilpotency class 3 or more attaining the bound \cite{Hatui}. She also gave an example of 3-group such that the bound is attained. 

In the following theorem we give a more direct proof of Hatui's result and classify 3-groups such that bound (\ref{eq0}) is attained completing the classification of such groups.

\begin{theorem}\label{thm2}
Let $G$ be a finite $p$-group of order $p^n$ with $|\gamma_2(G)| = p^k$. Then $|M(G)| = p^{\frac{1}{2}(n-k-1)(n+k-2)+1}$ if and only if $G$ is one of the following group.
\begin{itemize}
 \item[1.] $G_1 = E_p \times C_p^{(n-3)}$, where $E_p$ is the extraspecial $p$-group of order $p^3$ and exponent $p$ for an odd prime $p$,
  \item[2.] $G_2 = C_p^{(4)} \rtimes C_p$ for an odd prime $p$,
   \begin{align*} 3. \ G_3  = \Big\langle\alpha_1, \beta_1, \alpha_2, \beta_2, \alpha_3, \beta_3 \ | \ [\alpha_1, \alpha_2] = \beta_3, [\alpha_2, \alpha_3] = \beta_1, [\alpha_3, \alpha_1] = \beta_2, \\ [\alpha_i, \beta_j] = 1, \alpha_i^p = \beta_i^p = 1 \ (p \ \text{an odd prime})\ (i, j = 1,2,3)\Big\rangle.
\end{align*}
\begin{align*} 4. \ G_4  = \Big\langle\alpha_1, \beta_1, \alpha_2, \beta_2, \alpha_3, \beta_3, \gamma \ | \ [\alpha_1, \alpha_2] = \beta_3, [\alpha_2, \alpha_3] = \beta_1, [\alpha_3, \alpha_1] = \beta_2, \\ [\beta_i, \alpha_i] = \gamma, \ \ [\alpha_i, \beta_j] = 1, \alpha_i^3 = \beta_i^3 = 1 \ (i= 1,2,3, \ \ j = 2,3)\Big\rangle.
\end{align*}
\end{itemize}
\end{theorem}

\vspace{.6cm}

A finite $p$-group of order $p^n$ is said to be of maximal class if its nilpotency class is $n-1$. Let $G$ be finite $p$-group of maximal class and order $p^n$. Since $G$ is generated by 2 elements, it follows by a result of Gasch\"{u}tz \cite{Gaschutz} that $|M(G)| \leq p^{n-1}$. Moravec proved for $n > p+1$, that $|M(G)| \leq p^{\frac{p+1}{2}\left \lceil{\frac{n-1}{p-1}}\right \rceil}$ \cite{Moravec}. Improving his result we prove the following theorem.

\begin{theorem} \label{thm3}
Let $G$ be a finite $p$-group of maximal class and of order $p^n$ for an odd prime $p$ and $n \geq 4$.  Then $|M(G)| \leq p^{\frac{n}{2}}$.  
\end{theorem}

\section{Prerequisites}\label{sec2}
Let $G$ be a finte $p$-group of nilpotency class $c$ and $\overline{G}$ be the factor group $G/Z(G)$.  Define homomorphism
\[\Psi_2: \overline{G}^{ab} \otimes \overline{G}^{ab} \otimes \overline{G}^{ab} \mapsto \frac{\gamma_2(G)}{\gamma_{3}(G)} \otimes \overline{G}^{ab}\]
\[\text{by} \ \Psi(\overline{x_1} \otimes \overline{x_2} \otimes \overline{x_{3}}) 
= \overline{[x_1, x_2]} \otimes \overline{x_{3}} + \overline{[x_2, x_{3}]} \otimes \overline{x_1} + \overline{[x_3, x_{1}]} \otimes \overline{x_2}.\]

For $3 \leq i \leq c$ define homomorphisms 
\[\Psi_i: \underbrace{\overline{G}^{ab} \otimes \overline{G}^{ab} \cdots \otimes \overline{G}^{ab}}_{i+1 \ \text{times}} \ \ \mapsto \frac{\gamma_i(G)}{\gamma_{i+1}(G)} \otimes \overline{G}^{ab}\]
by
\begin{eqnarray*}
\Psi_i(\overline{x_1} \otimes \overline{x_2} \otimes \cdots \otimes \overline{x_{i+1}}) 
& = & \overline{[x_1, x_2, \cdots, x_i]_l} \otimes \overline{x_{i+1}} + \overline{[x_{i+1}, [x_1, x_2, \cdots x_{i-1}]_l]} \otimes \overline{x_i} \\
&& +\overline{[[x_i, x_{i+1}]_r, [x_1, \cdots, x_{i-2}]_l]} \otimes \overline{x_{i-1}} \\
&& + \overline{[[x_{i-1}, x_i, x_{i+1}]_r, [x_1, x_2, \cdots, x_{i-3}]_l]} \otimes \overline{x_{i-2}} \\
&& + \cdots + \overline{[x_2, \cdots, x_{i+1}]_r} \otimes \overline{x_1}\\
\end{eqnarray*}
where 
\[[x_1, x_2, \cdots x_i]_r = [x_1, [\cdots [x_{i-2},[x_{i-1},x_i]]\ldots]\]
and 
\[[x_1, x_2, \cdots x_i]_l = [\ldots[[x_1, x_2], x_3], \cdots, x_i].\]

\vspace{.2cm}
The following proposition was given by Ellis and Weigold \cite[Proposition 1 and the comments on page 192 following the proof of Theorem 2]{Ellis1}.

\begin{proposition}\label{prop1}
Let $G$ be a finite $p$-group and $\overline{G}$ be the factor group $G/Z(G)$. Then 
\[\Big{|}M(G)\Big{|}\Big{|}\gamma_2(G)\Big{|}\prod_{i=2}^{c}\Big{|}Im \Psi_i\Big{|} \leq \Big{|}M(G^{ab})\Big{|}\prod_{i=2}^{c}\Big{|}\frac{\gamma_i(G)}{\gamma_{i+1}(G)} \otimes \overline{G}^{ab}\Big{|}.\]
\end{proposition}

The following Lemma is from \cite{Rai1}.%gives a bound for the schur multiplier of an abelian $p$-group. 

\begin{lemma}\cite[Lemma 2.1]{Rai1}\label{lem1}
Let $G$ be an abelian $p$-group of order $p^n$ such that $G = C_{p^{\alpha_1}} \times C_{p^{\alpha_2}} \times \cdots \times C_{p^{\alpha_d}} (\alpha_1 \geq \alpha_2 \geq \cdots \geq \alpha_d)$ and $|G| = p^n$, then $|M(G)| \leq p^{\frac{1}{2}(d(G)-1)(n-(\alpha_1-\alpha_d))}$.
\end{lemma}

The following Lemma is from \cite{Khukhro}.

\begin{lemma}\cite[Lemma 3.6 (c)]{Khukhro}\label{lem2}
Let $G$ be a group such that $G = \gen{M}$. Then $\gamma_k(G)$ is generated by simple commutators of weight $\geq k$ in the elements $m^{\pm}, m \in M$.
\end{lemma}

% In 1956 Green proved that $|M(G)| \leq p^{\frac{1}{2}n(n-1)}$ \cite{Green}. Since then, Green's bound has been reproved and generalized by many mathematicians. Wiegold, in 1965, gave a bound on $|\gamma_2(G)|$ in terms of $|G/Z(G)|$ and rederived the Green's bound using the existence of representation groups \cite{Wiegold1}. In 1967 Gasch\"{u}tz et al., sharpening Green's bound, proved in \cite{Gaschutz} that 
%\begin{equation*}
%|M(G)| \leq |M(G^{ab})||\gamma_2(G)|^{d(G/Z(G))-1}. \label{bnd_Gtz}
%\end{equation*}
% The bound of Gasch\"{u}tz et al. was further generalized by Vermani in 1969 \cite{Vermani1}. He obtained their result as a corollary of the bound 
 %\begin{equation*} 
 %|M(G)| \leq \bigg{|}M\bigg(\frac{G}{\gamma_c(G)}\bigg)\bigg{|}\bigg{|}Hom\bigg(\frac{G}{Z_{c-1}(G)}, \gamma_c(G)\bigg)\bigg{|}\bigg{/}|\gamma_c(G)|, \label{bnd_vrmni1}
 %\end{equation*}
%where $c$ is the nilpotency class of $G$. This bound was reproved by Jones using a different method in \cite{Jones1}. In 1969 Green's bound was generalized by Wiegold  \cite{Wiegold2} when he proved that 

%This was reproved by Niroomand and Russo \cite{Niroomand1} using a different method. 

 %Moreover, They claimed to improve the bound when $G^{ab} \neq C_{p^e} \times C_{p^e} \times \ldots  \times C_{p^e}$ proving that 
%\begin{equation}
%|M(G)| \leq  p^{\frac{1}{2}(d_G-1)(n+k-1)} \label{bnd_nrmnd1} 
%\end{equation}
%in this case. 

\section{Proofs of Theorems}

{\bf Proof of Theorem \ref{thm1}}
Let $\Psi_i$ be as defined in Section \ref{sec2} and $d(G/Z(G)) = \delta$. Following Proposition \ref{prop1} we have that 
\[|M(G)||\gamma_2(G)|\prod_{i=2}^{c}|Im \Psi_i| \leq |M(G^{ab})|p^{k\delta}.\]
Applying Lemma \ref{lem1} this gives
\[|M(G)|\prod_{i=2}^{c}|Im \Psi_i| \leq |p^{\frac{1}{2}(d-1)(n-k)+k(\delta-1)}.\]
so that
\begin{equation}\label{eq1}
|M(G)|\prod_{i=2}^{c}|Im \Psi_i| \leq |p^{\frac{1}{2}(d-1)(n+k) - k(d-\delta)}.
\end{equation}
Choose a subset $S= \{x_1, x_2, \ldots, x_{\delta}\}$ of $G$ such that $\{\overline{x_1},\overline{x_2}, \cdots, \overline{x_{\delta}}\}$ be a minimal generating set for $G/Z(G)$. Fix $i \leq \min(\delta, c)$. Since $ i \leq c$, $\gamma_i(G)/\gamma_{i+1}(G)$ is a non-trivial group. Using Lemma \ref{lem2} we can choose a commutator $[y_1,y_2, \cdots, y_i]$ of weight $i$ such that $[y_1,y_2, \cdots, y_i] \notin \gamma_{i+1}(G)$ and $y_1, \ldots, y_i \in S$ . Since $i \leq \delta$, $S \backslash \{y_1,y_2, \cdots, y_i\}$ contains at least $\delta - i$ elements. Choose any $\delta - i$ elements $z_1, z_2, \ldots, z_{\delta - i}$ from $S \backslash \{y_1,y_2, \cdots, y_i\}$. Since $[y_1,y_2, \cdots, y_i] \notin \gamma_{i+1}(G)$ and $z_j \notin \{y_1,y_2, \cdots, y_i\}$, $\Psi_i(\overline{y_1}, \ldots, \overline{y_i}, \overline{z_j}) \neq 1$. Notice that the set $\{\Psi_i(\overline{y_1}, \ldots, \overline{y_i}, \overline{z_j}) \ \ | \ \ 1 \leq j \leq \delta -i\}$ is a minimal generating set for $\gen{\{\Psi_i(\overline{y_1}, \ldots, \overline{y_i}, \overline{z_j}) \ \ | \ \ 1 \leq j \leq \delta -i\}}$ because $\{\overline{x_1},\overline{x_2}, \cdots, \overline{x_{\delta}}\}$ is a minimal generating set for $G/Z(G)$. It follows that $|\text{Im}\Psi_i| \geq p^{\delta - i}$. Putting this in Equation \ref{eq1} we get the required result.

 \vspace{.3cm}
     
{\bf Proof of Theorem \ref{thm2}}
Let $|M(G)| = p^{\frac{1}{2}(n-k-1)(n+k-2)+1}$. In view of \cite[Theorem 1.1]{Rai2} suppose that the nilpotency class of $G$ is at least 3. From the following exact sequence  \cite[Corollary 3.2.4 (ii)]{Karpilowski}
\[1 \mapsto X \mapsto G/\gamma_2(G) \otimes \gamma_c(G) \mapsto M(G) \mapsto M(G/\gamma_c(G)) \mapsto \gamma_c(G) \mapsto 1,\]
it follows that the Bound \ref{eq0} is attained for the group $G/\gamma_c(G)$. Applying induction the bound \ref{eq0} is attained for $G/\gamma_3(G)$. But $G/\gamma_3(G)$ is of nilpotency class 2, therefore by \cite[Theorem 1.1]{Rai2}, $p \neq 2$. 
Also, it follows from Theorem \ref{thm1} that $d(G) \leq 3$. Now using \cite[Theorem 2.2]{Niroomand2} we see that $d(G) = 3$. Let $\Psi_i$ be the maps as given in Section \ref{sec2}. By simplyfying notations 
\begin{eqnarray*}
\Psi_3(\overline{x_1} \otimes \overline{x_2} \otimes \overline{x_3} \otimes \overline{x_{4}}) 
& = & \overline{[[x_1, x_2], x_3]} \otimes \overline{x_{4}} + \overline{[x_{4}, [x_1, x_2]]} \otimes \overline{x_3} 
+\overline{[[x_3, x_{4}], x_1]} \otimes \\ 
&& \overline{x_{2}} + \overline{[x_2, [x_3, x_{4}]]} \otimes \overline{x_{1}}. \\
\end{eqnarray*}

As in the proof of Theorem \ref{thm1} we have $|\text{Im} \Psi_2| \geq p^{\delta -2}$. 
Applying Equation \ref{eq1} we get that $d = \delta$. It follows from Equation \ref{eq1} that 
\begin{equation}\label{eq3}
|M(G)|\prod_{i=3}^{c}|Im \Psi_i| \leq p^{\frac{1}{2}(d-1)(n+k-2)+1}.
\end{equation}

Therefore $\text{Im} \Psi_3 = \{1\}.$ Let $G$ be generated by $\alpha_1, \alpha_2, \alpha_3$. Then for $i \neq j$ 
\[\Psi_3(\overline{\alpha_i} \otimes \overline{\alpha_j} \otimes \overline{\alpha_i} \otimes \overline{\alpha_j}) = 2(\overline{[\alpha_i, \alpha_j, \alpha_i]} \otimes \overline{\alpha_j}) + 2(\overline{[\alpha_j, [\alpha_i, \alpha_j]]} \otimes \overline{\alpha_i)}.\]
This shows that $[\alpha_i, \alpha_j, \alpha_i] \in \gamma_4(G)$ because $p \neq 2$.
Now for $i \neq j \neq k \neq i$, consider 
\[\Psi_3(\overline{\alpha_i} \otimes \overline{\alpha_j} \otimes \overline{\alpha_k} \otimes \overline{\alpha_i}) = \overline{[\alpha_i, \alpha_j, \alpha_k]} \otimes \overline{\alpha_i} + \overline{[\alpha_j, [\alpha_k, \alpha_i]]} \otimes \overline{\alpha_i}.\]
Therefore
\[[\alpha_i, \alpha_j, \alpha_k]\gamma_4(G) = [\alpha_k, \alpha_i, \alpha_j]\gamma_4(G).\]
Putting $(i, j, k) = (1,2,3)$ and $(2,3,1)$ gives 
\[[\alpha_1, \alpha_2, \alpha_3]\gamma_4(G) = [\alpha_3, \alpha_1, \alpha_2]\gamma_4(G).\]
and
\[[\alpha_2, \alpha_3, \alpha_1]\gamma_4(G) = [\alpha_1, \alpha_2, \alpha_3]\gamma_4(G).\]
respecticely.

Applying Hall-Witt identity we see that $[\alpha_2, \alpha_3, \alpha_1]^3 \in \gamma_4(G)$. Since $[\alpha_i, \alpha_j, \alpha_i] \in \gamma_4(G)$ and $\gamma_3(G)/\gamma_4(G)$ is non-trivial, we have $[\alpha_1, \alpha_2, \alpha_3] \notin \gamma_4(G)$. It follows that $p =3$. Let $G$ be a group of nilpotency class at least 4. 
Consider the map $\Psi_4$. By simplyfying notations 

\begin{eqnarray*}
\Psi_4(\overline{x_1} \otimes \overline{x_2} \otimes \overline{x_3} \otimes \overline{x_{4}} \otimes \overline{x_5}) & = & \overline{[x_1, x_2, x_3, x_4]} \otimes \overline{x_{5}} + \overline{[x_{5}, [x_1, x_2, x_3]]} \otimes \overline{x_4}  \\
&& +\overline{[[x_4, x_{5}], [x_1, x_2]]} \otimes \overline{x_{3}} + \overline{[[x_3, [x_4, x_5]], x_1]} \otimes \overline{x_2}. \\
&& + \overline{[x_2, [x_3, [x_4, x_5]]]} \otimes \overline{x_1}. \\
\end{eqnarray*}

Since $\gamma_4(G)/\gamma_5(G)$ is non-trivial, one of the elements $[\alpha_1, \alpha_2, \alpha_3, \alpha_i], \ \ i = 1,2,3$ does not belong to $\gamma_5(G)$. Suppose $[\alpha_1, \alpha_2, \alpha_3, \alpha_1] \notin \gamma_5(G)$. Then $\Psi_4(\overline{\alpha_1} \otimes \overline{\alpha_2} \otimes \overline{\alpha_1} \otimes \overline{\alpha_2} \otimes \overline{\alpha_3})$ is non-identity so that Im$\Psi_4$ is non-trivial. 
Similarly supposing $[\alpha_1, \alpha_2, \alpha_3, \alpha_2] \notin \gamma_5(G)$, the element $\Psi_4(\overline{\alpha_1} \otimes \overline{\alpha_2} \otimes \overline{\alpha_2} \otimes \overline{\alpha_1} \otimes \overline{\alpha_3})$, while supposing $[\alpha_1, \alpha_2, \alpha_3, \alpha_3] \notin \gamma_5(G)$, the element $\Psi_4(\overline{\alpha_1} \otimes \overline{\alpha_2} \otimes \overline{\alpha_3} \otimes \overline{\alpha_3} \otimes \overline{\alpha_1})$ give that Im$\Psi_4$ is non-trivial. This, in view of Equation \ref{eq3}, gives a contradiction. Therefore $G$ is a 3-group of nilpotency class 3. Hence we have  
\[[\alpha_1, \alpha_2, \alpha_3] = [\alpha_3, \alpha_1, \alpha_2] = [\alpha_2, \alpha_3, \alpha_1].\]
Since $[\alpha_1, \alpha_2, \alpha_3] \neq 1$ we get that $[\alpha_i, \alpha_j] \notin \gamma_3(G)$ for $i, j = 1,2,3$. Also, since $[\alpha_i, \alpha_j, \alpha_i] = 1$, it follows that $[\alpha_i, \alpha_j]\gamma_3(G)$ can not be generated by $\{[\alpha_k, \alpha_l]\gamma_3(G) \ \ | \\ \ \ (k,l) \neq (i,j) \ \ \text{or} \ \ (j,i)\}$. This shows that $\gamma_2(G)/\gamma_3(G)$ is generated by 3 elements. Using Equation \ref{eq3} $G/\gamma_2(G)$ is elementary abelian. So that $\gamma_2(G)/\gamma_3(G)$ is elementary abelian. Hence $|\gamma_2(G)/\gamma_3(G)| = 3^3$. Therefore $|G| = 3^7$. Now it can be checked using GAP that the bound is attained if and only if $G = G_4$. This completes the proof. %no group of order $3^7$ with $|\gamma_2(G)| = 3^4$ and exponent more than 3 attains the bound. Therefore $G$ is given by the following presentation

%\begin{align*} \ G_4  = \Big\langle\alpha_1, \beta_1, \alpha_2, \beta_2, \alpha_3, \beta_3, \gamma \ | \ [\alpha_1, \alpha_2] = \beta_3, [\alpha_2, \alpha_3] = \beta_1, [\alpha_3, \alpha_1] = \beta_2, \\ [\beta_i, \alpha_i] = \gamma, \ \ [\alpha_i, \beta_j] = 1, \alpha_i^3 = \beta_i^3 = 1 \ \ (i= 1,2,3, \ \ j = 2,3)\Big\rangle.
%\end{align*}

%Conversely, it can can be checked using GAP that the bound \ref{eq0} is attained for the group $G_4$.

\vspace{.3cm}

{\bf Proof of Theorem \ref{thm3}} Let $P_1 = C_G(\gamma_2(G)/\gamma_4(G)).$ Choose arbitrary elements $s \in G \backslash P_1 \cup C_G(\gamma_{n-2}(G))$ and $s_1 \in P_1 \backslash \gamma_2(G)$. Then $s$ and $s_1$ generate $G$. If we define $s_i = [s_{i-1}, s]$ for $i \geq 2$, then $s_i \in \gamma_i(G) \backslash \gamma_{i+1}(G)$. Let $\Psi_i, i \geq 3$, be the map as defined in Section \ref{sec2}. Then
\begin{eqnarray*}
\Psi_i(\overline{s_1} \otimes \overline{s} \otimes \overline{s} \otimes \cdots \otimes \overline{s} \otimes \overline{s_1}) & = & \overline{[s_1, s, s, \cdots,s]_l}\overline{[s, s, \cdots, s, s_1]_r} \otimes \overline{s_1} + \overline{t} \otimes \overline{s}
\end{eqnarray*}
for some $t \in G$. \\
Notice, for an odd $i$, that 
\[\overline{[s_1, s, s, \cdots,s]_l} = \overline{[s, s, \cdots, s, s_1]_r}.\]
Since $p \neq 2$, it follows that $\Psi_i(\overline{s_1} \otimes \overline{s} \otimes \overline{s} \otimes \cdots \otimes \overline{s} \otimes \overline{s_1})$ is non-identity so that Im$\Psi_i$ is non-trivial. Using this fact Equation \ref{eq1} gives the required result.

%We use the following notations. For a multiplicatively written group $G$ let $x,y,a \in G$. Then $[x,y]$ denotes the commutator $x^{-1}y^{-1}xy$ and $x^a$ denotes the conjugate of $x$ by $a$ i.e. $a^{-1}xa$. 
%By $\gen{x}$ we denote the cyclic subgroup generated by $x$. %By $Z(G)$ we denote the center of $G$. % and by $C_G(H)$ we denote the centralizer of $H$ in $G$ where $H$ is a %subgroup of $G$. 
%We write $\gamma_2(G)$ for the commutator subgroup of $G$. %The group of all homomorphisms from a group $H$ to an abelian group $K$ is denoted by $Hom(H, K)$. 
%For $x \in G$,
 % %\ x^G$ denotes the $G$-conjugacy class of $x$ and 
%$[x, G]$ denotes the set of all $[x,y]$, for $y \in G$. %By $C_p$ we denote the cyclic group of order $p$. 
%The subgroup generated by all the elements of order $p$ is denoted by $\Omega_1(Z(G)$. By a class $c$ group we mean a group of nilpotency class $c$.  \\

\vspace{.3cm}

{\bf Acknowledgements:} I am very grateful to my post-doctoral superviser Prof. Boris Kunyavski\u{\i} for his encouragement and support. This research was
supported by Israel Council for Higher Education’s fellowship program and by ISF
grant 1623/16.

\end{document}